\documentclass[leqno]{amsart}

\usepackage{amssymb,latexsym}
\usepackage{amsmath}
\usepackage{graphicx} 
\usepackage{enumerate}
\usepackage{newlattice}

\theoremstyle{plain}

\DeclareMathOperator{\threeC}{3C}

\newcommand{\col}[1]{\tup{col}(#1)}

\newcommand{\swing}{\mathbin{\raisebox{2.0pt}
       {\rotatebox{160}{$\curvearrowleft$}}}}
\newcommand{\exswing}{\stackrel{\textup{ ex}}{\swing}}
\newcommand{\inswing}{\stackrel{\textup{\,in}}{\swing}}

\newcommand{\Kpeek}{K^{\includegraphics[scale = .2]{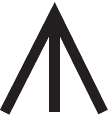}}}

\theoremstyle{definition}

\newcommand{\seven}{\SfS 7}

\newcommand{\mthree}{\msf{M}_3}

\begin{document}
\title[On a property of congruence lattices]{On a property of congruence lattices
of slim, planar, semimodular lattices}
\author[G.\ Gr\"atzer]{George Gr\"atzer}
\email{gratzer@mac.com}
\urladdr{http://server.maths.umanitoba.ca/homepages/gratzer/}
\address{University of Manitoba}
\date{\today}

\begin{abstract} 
In a 2121 paper with G\'abor Cz\'edli,
we introduced and verified the Three-pendant Three-crown Property, 3P3C,
for congruence lattices of slim (no $\mthree$ sublattice), planar, semimodular lattices.
The proof is very long; 
in part, because it relies on Cz\'edli's 2021 paper on lamps.

This paper verifies 3P3C using the Swing Lemma, 
an elementary and short approach.
\end{abstract}

\maketitle

\section{Introduction}\label{S:Introduction}

\begin{figure}[htb]
\centerline{\includegraphics{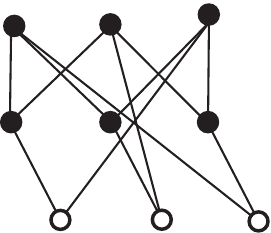}}
\caption{The ordered set Three-pendant Three-crown, $R_3$;\\ the elements of the crown are black-filled}
\label{F:3P3C}
\end{figure}

In my joint paper with G. Cz\'edli~\cite{CG21}, 
we defined the \emph{Three-pendant Three-crown Property} (\emph{3P3C Property}),
as follows:
\begin{enumeratei}
\item[] \emph{The ordered set $R_3$ of Figure~\ref{F:3P3C} has no cover-preserving
embedding into the ordered set of join-irreducible congruences of $K$, $\Ji{\Con K}$.}
\end{enumeratei}

The following is the result of the same paper.

\begin{named}{3P3C Theorem}
Let $K$ be a slim, planar, semimodular lattice. 
Then $\Con K$ satisfies the 3P3C Property. 
\end{named}

In this paper, we provide a short and elementary proof, 
utilizing the Swing Lemma of my paper~\cite{gG15}.

This paper is largely self-contained. 
Apart from some very elementary concepts 
(semimodular lattice, complemented lattice, interval, etc.),
Sections~\ref{S:basic} and \ref{S:Swing} present all the concepts and results we need,
sometimes informally.

For more extensive background, download 
Part~I of the book \cite{CFL3}, see 

\smallskip

{\tt arXiv:2104.06539}

\subsection*{Declarations}
Data sharing is not applicable to this article as no datasets were generated
or analyzed during the current study. 
I further state that there is no conflict of interest.
\subsection*{Funding statement} There was no funding for this project.

\section{Some basic concepts and results}\label{S:basic}

\subsection{Join-irreducible congruences}\label{S:xx}

Let $L$ be a finite lattice. 
In $L$, let us call an interval $[a, b]$  \emph{prime}, 
if  $[a, b]$ has only two elements, namely, $a$ and~ $b$. 
A prime interval $[a, b]$ of $L$ is represented by an \emph{edge} $E$
in a diagram of $L$. We will use the terms prime interval and edge interchangeably.

If $[a, b]$ is a prime interval, 
then $\con{a,b}$, the smallest congruence collapsing it, 
is a join-irreducible congruence, and conversely. Similarly, for edges.

\subsection{Two acronyms}\label{S:acronyms}

Let $L$ be a planar lattice. 
The lattice $L$ is \emph{slim}, if it has no~$\mthree$ sub\-lattice.
The acronym SPS stands for slim, planar, semimodular, as in SPS lattices.

Following my joint paper with E.~Knapp \cite{GKn07},
a planar semimodular lattice~$L$ is \emph{rectangular}, if its left boundary chain
has exactly one doubly-irreducible element,~$c_l$,
and its right boundary chain 
has exactly one doubly-irreducible element, $c_r$,
and these elements are complementary, that~is,
\begin{align*}
c_l \jj c_r &=1,\\
c_l \mm c_r &=0.
\end{align*}

The acronym SR stands for slim rectangular, as in SR lattices.

\subsection{Forks}\label{S:Forks}

\begin{figure}[b!]
\centerline{\includegraphics[scale=1]{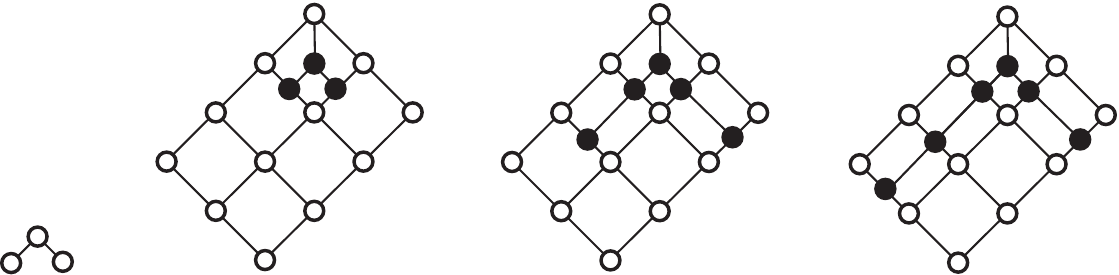}}
\caption{Inserting a fork}\label{F:newfork}
\end{figure}

Let $L$  be an SPS lattice. 
Following G.~Cz\'edli and E.\,T. Schmidt~\cite{CS12},
we introduce the concept of \emph{inserting a fork into $L$ at $C$}, 
where $C$ is an interval of~$L$ that is a covering square.
We start with the ordered set~$F$, the fork, as pictured in the first diagram of Figure~\ref{F:newfork} and the SPS lattice~$L = \SC 3 \times \SC 4$, 
with the top interval of $L$ as $C$.
We place $F$ into $C$ as in the second diagram of Figure~\ref{F:newfork}. 
We~obtain a lattice but it is not semidistributive, 
so the left bottom element of the fork requires an additional element,
left and down from it, to correct, and symmetrically.
The~lattice we obtain, as in the third diagram of Figure~\ref{F:newfork}, 
is still not semidistributive, so the left bottom element of the new fork
similarly requires an additional element to correct.

It is easy to see that we obtain an SPS lattice. 

\subsection{Representation theorems}\label{S:xx}

In this branch of lattice theory, a \emph{representation theorem} (abbreviated as RT)
is a statement of the form: every finite distributive lattice
can be represented as the congruence lattice of a lattice of some type.
The~first RT found is the following.

\begin{named}{Basic Representation Theorem} 
Every finite distributive lattice $D$ can be
represented as the congruence lattice of a finite lattice $L$.
\end{named}

This is due to R.\,P. Dilworth from around 1944 (see the book \cite{BFK90}).
The Basic RT was not published until 1962
in my joint paper with E.\,T. Schmidt \cite{GS62}.

\begin{figure}[htb]
\centerline{\includegraphics{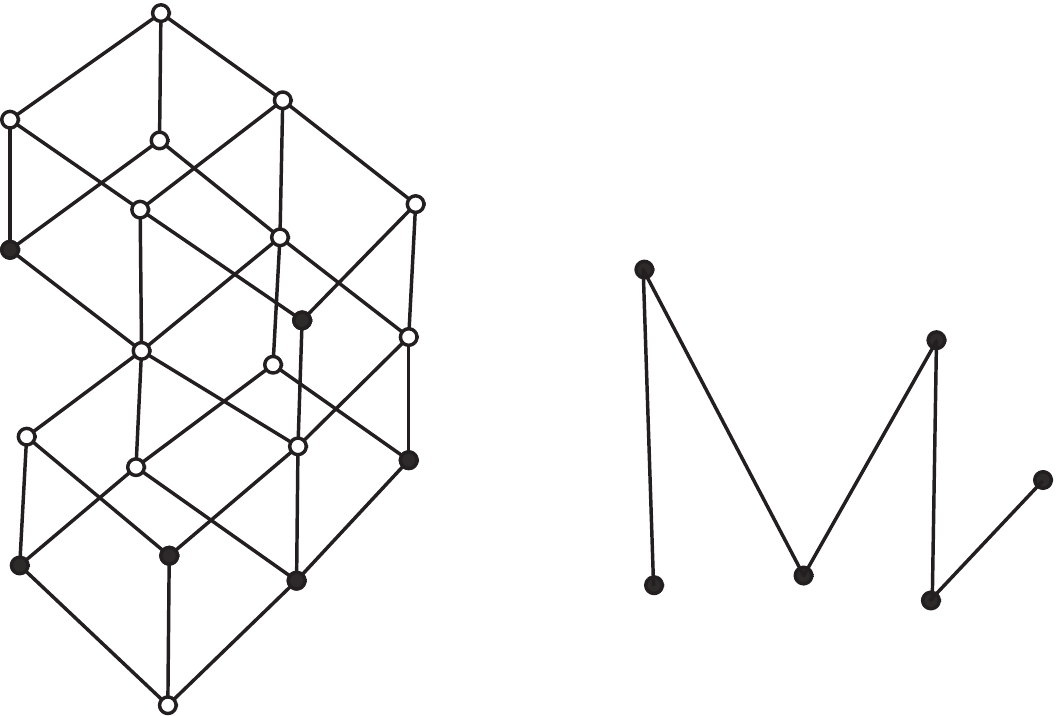}}
\caption{A finite distributive lattice $D$ 
and the ordered set of join-irreducible elements $P$ of $D$}
\label{F:free-new}
\end{figure}

With a finite distributive lattice $D$, we can associate the ordered set $P$
of join-irreducible elements; also, from the ordered set $P$, we can easily 
reconstruct $D$. This may reduce a complex finite distributive lattice $D$
to a much smaller ordered set $P$ of simple structure, 
see Figure~\ref{F:free-new} for an illustration.

So we may rephrase the Basic Representation Theorem as follows.

\begin{named}{Basic Representation Theorem'} 
Every finite ordered set $P$ can be
represented as the ordered set, $\Ji {\Con L}$, of join-irreducible congruences 
of a finite lattice~$L$.
\end{named}

The first specialized RT was in the same paper.
Recall that a finite lattice $L$ is \emph{sectionally complemented}, if 
every ideal of $L$ is a complemented lattice, that is, for all $b \leq a \in L$,
there exists a $c \in L$ such that $a = b \jj c$ and $0 = b \mm c$.

\begin{named}{RT for sectionally complemented latticesm} 
Every finite distributive lattice~$D$ can be
represented as the congruence lattice of a finite sectionally complemented lattice $L$
\end{named}

\subsection{Two-cover Theorem}\label{S:Two}

The topic we make a contribution to in this paper 
started with the following result in my paper \cite{gG16}.
Let $K$ be an SPS lattice with at least three elements.
A finite ordered set~$P$ satisfies the 
 \emph{Two-Cover Condition}\index{Two-Cover Condition},
if any element of~$P$ has at most two covers.

\begin{named}{Two-cover Theorem}
The ordered set of join-irreducible congruences 
of an SPS lattice $K$ satisfies the Two-Cover Condition.
\end{named}

\section{The Swing Lemma}\label{S:Swing}

\begin{figure}[htb]
\centerline{\includegraphics{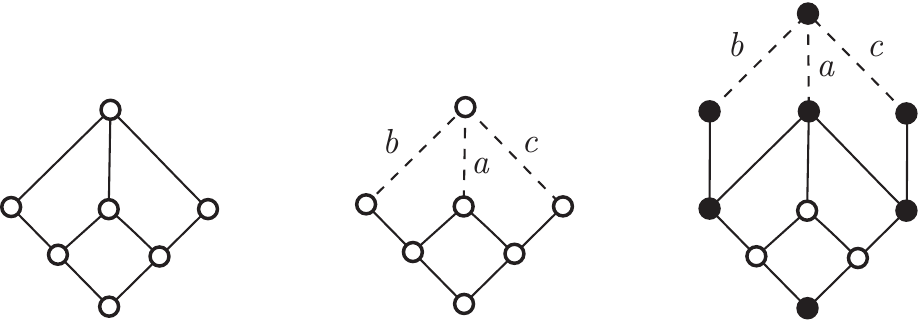}}
\caption{The $\seven$ lattice and two peak sublattices $S(a, b, c)$; 
in~the second example, the elements of the peak sublattice are black filled}
\label{F:V}
\end{figure}

The crucial tool employed in this paper is the Swing Lemma, introduced in this section.
Let $K$ be an SPS lattice. 
The $\seven$ lattice is illustrated in Figure~\ref{F:V}.
A~\emph{peak~$\seven$ sublattice} (\emph{peak sublattice}, for short)
$S$ is an $\seven$ sublattice,
whose top three edges are covering. See Figure~\ref{F:V} for examples.

\begin{figure}[htb!]
\centerline{\includegraphics{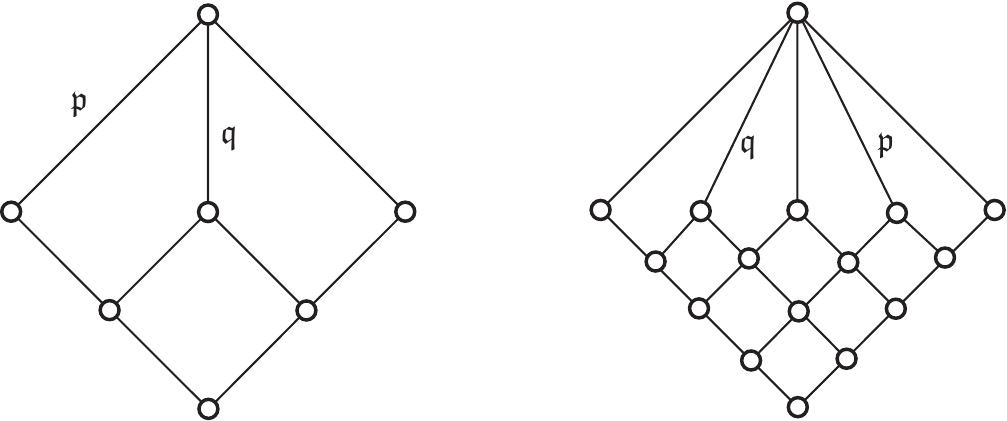}} 
\caption{Swings, $\fp \protect\swing \fq$}\label{F:Swings}
\end{figure}

For the prime intervals $\fp, \fq$ of a slim, planar, semimodular lattice (SPS lattice)~$L$,
we define a binary relation:
$\fp$~\emph{swings to} $\fq$, written as $\fp \swing \fq$,
if $1_\fp = 1_\fq$, 
this element covers at least three elements,
and $0_\fq$ is neither the left-most nor the right-most element
covered by $1_\fp = 1_\fq$, see Figure~\ref{F:Swings} for two examples. 
If $0_\fp$ is either the left-most or the right-most element
covered by $1_\fp = 1_\fq$, then we call the swing \emph{external},
in formula, $\fp \exswing \fq$. 
Otherwise, the swing is \emph{internal},
in formula, $\fp \inswing \fq$.
Figure~\ref{F:Swings} shows an external swing, $\fp \exswing \fq$; 
and an internal swing, $\fp \inswing \fq$.

This paper is a series of applications of the following result (see my paper \cite{gG15}).

\begin{named}{Swing Lemma}
Let $L$ be an SPS lattice 
and let $\fp$ and $\fq$ be distinct prime intervals in $L$. 
Then $\fq$ is collapsed by $\con{\fp}$ if{}f 
there exists a prime interval~$\fr$ 
and a sequence of pairwise distinct prime intervals
\begin{equation}\label{Eq:sequence}
\fr = \fr_0, \fr_1, \dots, \fr_n = \fq
\end{equation}
such that $\fp$ is up perspective to $\fr$ and 
$\fr_i$ is down perspective to or swings to $\fr_{i+1}$
for $i = 0, \dots, n-1$. 
In addition, the sequence \eqref{Eq:sequence} is descending\cou  
\begin{equation}\label{E:geq}
   1_{\fr_0} \geq 1_{\fr_1} \geq \dots \geq 1_{\fr_n}.
\end{equation}
\end{named}

Recall that $\fp$ is \emph{up perspective} to $\fr$,
if $1_\fp \jj 0_\fq = 1_\fq$ and $1_\fp \mm 0_\fq = 0_\fp$.
We define  \emph{down perspective} dually.

Most of the time, we do need need the full Swing Lemma, 
only the  following two special cases (which, combined, prove the Swing Lemma).

\begin{figure}[htb]
\centerline{\includegraphics{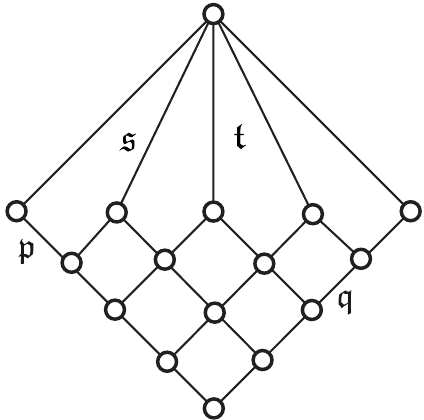}}
\caption{Illustrating the Equality Lemma, $\con{\fp} = \con\fq$}\label{F:equal}
\end{figure}

\begin{named}{Equality Lemma}
Let $\fp$ and $\fq$ be  prime intervals of $K$.
Then the equality \[\con\fp  = \con\fq\] holds if{}f  
there exist prime intervals $\fs$ and ${\ft}$ in~$K$, such that 
\begin{equation}\label{E:7(i)}
   {\fp}  \perspup \fs,\ \fs  \inswing \ft,\ \ft \perspdn \fq,
\end{equation}
see Figure~\ref{F:equal}.
\end{named}

\begin{named}{Covering Lemma} 
 Let $\fp$ and $\fq$ be  prime intervals of $K$.
The covering 
\[
   \con\fq  \prec \con{\fp}
\] 
holds in $\Ji{\Con K}$ if{}f there exist prime intervals $\fr,\ \fs,\ \ft, \ \fu,\  \fv$ in~$K$, such that 
\begin{equation}\label{E:8(i)}
   \fp  \perspup \fr,\  \fr \inswing \fs,\  \fs \perspdn \ft,\ 
      \ft \exswing \fu,\  \fu \perspdn \fq.
\end{equation}
In \eqref{E:8(i)}, we may have 
\begin{equation}\label{E:xx}
\fp  = \fr,\  \fr = \fs,\  \fs = \ft,\  \fu = \fq,
\end{equation}
or any combination thereof.
Moreover, there are prime intervals, 
$\ft,\fu$ and a peak sublattice $S(\ft,\fu)$ in $K$
with $\col \fu = \col \fq$ and $\col \ft = \col \fp$. 
\end{named}

\begin{figure}[htb]
\centerline{\includegraphics[scale=.16]{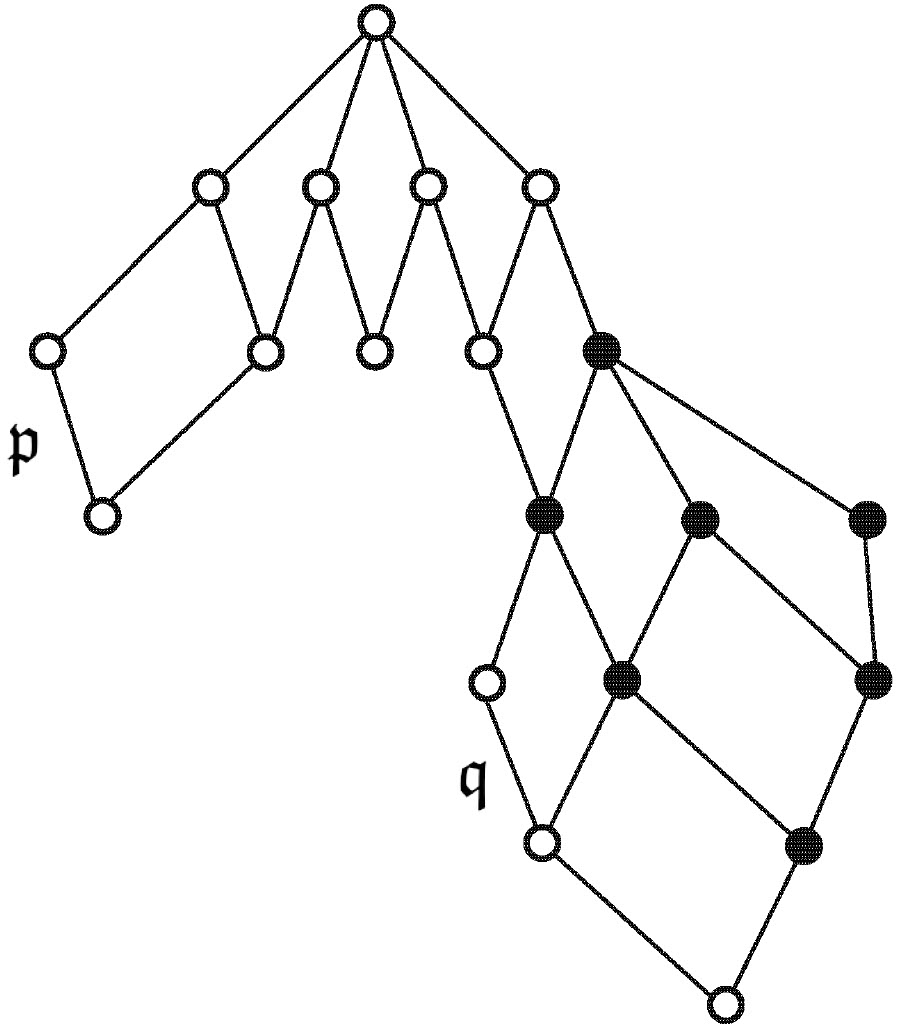}}
\caption{Illustrating the Covering Lemma, $\con\fq \prec \con{\fp}$;
the elements of the associated peak sublattice $S(\fp, \fq)$ are black filled}
\label{F:cover}
\end{figure}

For every  $\con\fq  \prec \con{\fp}$,
the Covering Lemma provides a peak sublattice.
This gives an approach to the proof of the 3P3C Theorem:
the covers in the ordered set~$R_3$ of Figure~\ref{F:3P3C} give us the peak sublattices
$S_1, S_2, \dots, S_{12}$. 
I~have not been able to use this approach;
$12$ peak sublattices are too many to handle. 
The key to the  proof is the
V-Lemma, which cuts the $12$ peak sublattices to $6$ (really, to $3$).

\section{The Problem}\label{S:Problem}

In my paper \cite{gG16}, I raised the the following.

\begin{named}{Problem}
Characterize congruence lattices of SPS lattices.
\end{named}

This paper is a contribution to it.

G. Cz\'edli maintains a list of related papers 
(mostly by Cz\'edli and myself), 56 as of this writing, see

\verb+http://www.math.u-szeged.hu/~czedli/m/listak/publ-psml.pdf+

\smallskip

\section{Preparing for the proof}\label{Proof}

The 3P3C theorem will be proved in Section~\ref{S:theorem}.
This section prepares for the proof by introducing four relations 
on the set of join-irreducible congruences of an SR lattice $K$.

\subsection{The V-relation}\label{S:Lemmas}

Let $K$ be an SR lattice.
The set of colors of $K$ is defined as the set of join-irreducible congruences of $K$,
that is, $\Ji{\Con K}$.
We define the ternary relation $V$
on the colors of $K$  as follows:
$V(a,b,c)$ holds in $K$, if $a, b, c$  are three \emph{distinct colors} and  $a \prec b, c$
(the first the smallest). We call this relation the \emph{V-relation}.
We denote by $\Kpeek$ those elements of $K$ 
that cover three or more elements.  

\begin{figure}[htb]
\centerline{\includegraphics{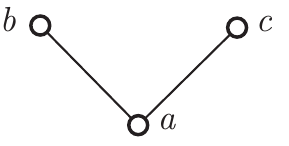}}
\caption{The V-relation on the colors of $K$, $\Ji{\Con K}$}
\label{F:Vrel}
\end{figure}

\newpage

\begin{named}{V-lemma}\hfill

\begin{enumeratei}
\item If $V(a,b,c)$ holds in $K$, then there is a peak sublattice $S(a, b, c)$, 
where $a$ colors the middle edge, 
$b, c$ color the other two top edges.
\item Let $S(a, b, c)$  be a peak sublattice of $K$,
as illustrated in Figure~\ref{F:V}, where $a$ colors the middle edge, 
$b, c$ color the other two top edges of $S(a, b, c)$.
Then $V(a,b,c)$ holds in $K$.
\item We make the following assumptions\cou
\begin{enumeratea}
\item  $V(a,b,c)$ holds in $K$\scu
\item $S(\fp, \fq)$ is a peak sublattice of $K$,
where $\fu$ is the middle edge, as in Figure~\ref{F:cover}\scu
\item $\col \fu = a$ and  $\col \ft = b$.
\end{enumeratea}

Then 
\[
   V(a,b,c) \mapsto 1_{S(\fp, \fq)}
\]
maps the set of V-relations of $K$ onto the subset $\Kpeek$ of $K$.
\end{enumeratei}
\end{named}

\begin{proof}
If $V(a,b,c)$ holds in $K$,
then $a \prec b$, so by the Covering Lemma, 
there are the prime intervals $\ft,\fu$ and the peak sublattice $S(\fp,\fq)$ in $K$
with $\col \fu = a$, $\col \ft = b$.
Let~$\ft'$ be the third top prime interval of $S(\fp,\fq)$ and let $b' = \col{\ft'}$.
Then $a \prec b, c,  b'$ and $b \neq c$.
By the Two-cover Theorem, $c = b'$. 
So the sublattice $S(\ft,\fu)$ is the peak sublattice  $S(\fp,\fq)$ we required.

(ii) and (iii) are trivial.
\end{proof}

\subsection{The W-relation}\label{S:W}

The \emph{W-relation}, $W$, in~$K$ 
is a $6$-ary relation on the set of colors of~$K$; 
the relation $V(a,b,c,d,e,f)$ holds in~$K$, if $V(a,b,c)$, $V(d,e,f)$, and $c = e$ hold in $K$.

Let $A$ and $B$ be edges in $K$
and let $Q$ denote the quadrilateral spanned by them. 
Let us say that $A$ \emph{faces} $B$, 
if they are opposite sides in $Q$.

\begin{figure}[htb]
\includegraphics[scale=.5]{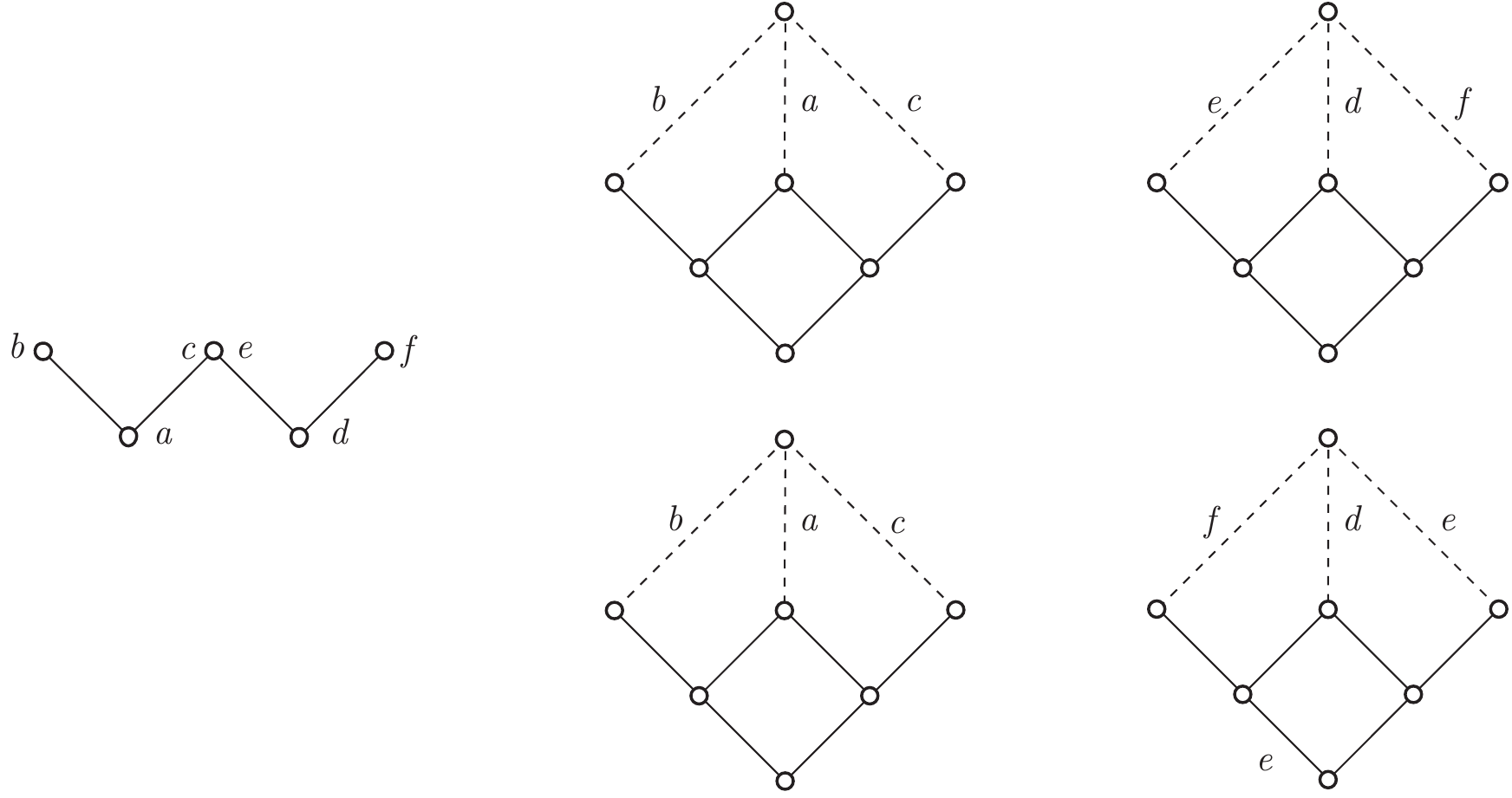}
\caption{Setting up the $W$-relation, two variants }\label{F:Wstart} 
\end{figure}

As shown in Figure~\ref{F:Wstart}, there are two possibilities 
how the two peak sublattices interrelate:

\begin{enumeratei}
\item Version 1. The edge colored by $c$ faces the edge colored by $e$.
\item Version 2. The edge colored by $c$ faces the edge colored by $f$.
\end{enumeratei}

The two variants are illustrated in Figure~\ref{F:W}.

\subsubsection{The W-relation, Version 1}\label{S:Wfirst}

\begin{figure}[htb]
\centerline{\includegraphics[scale=.6]{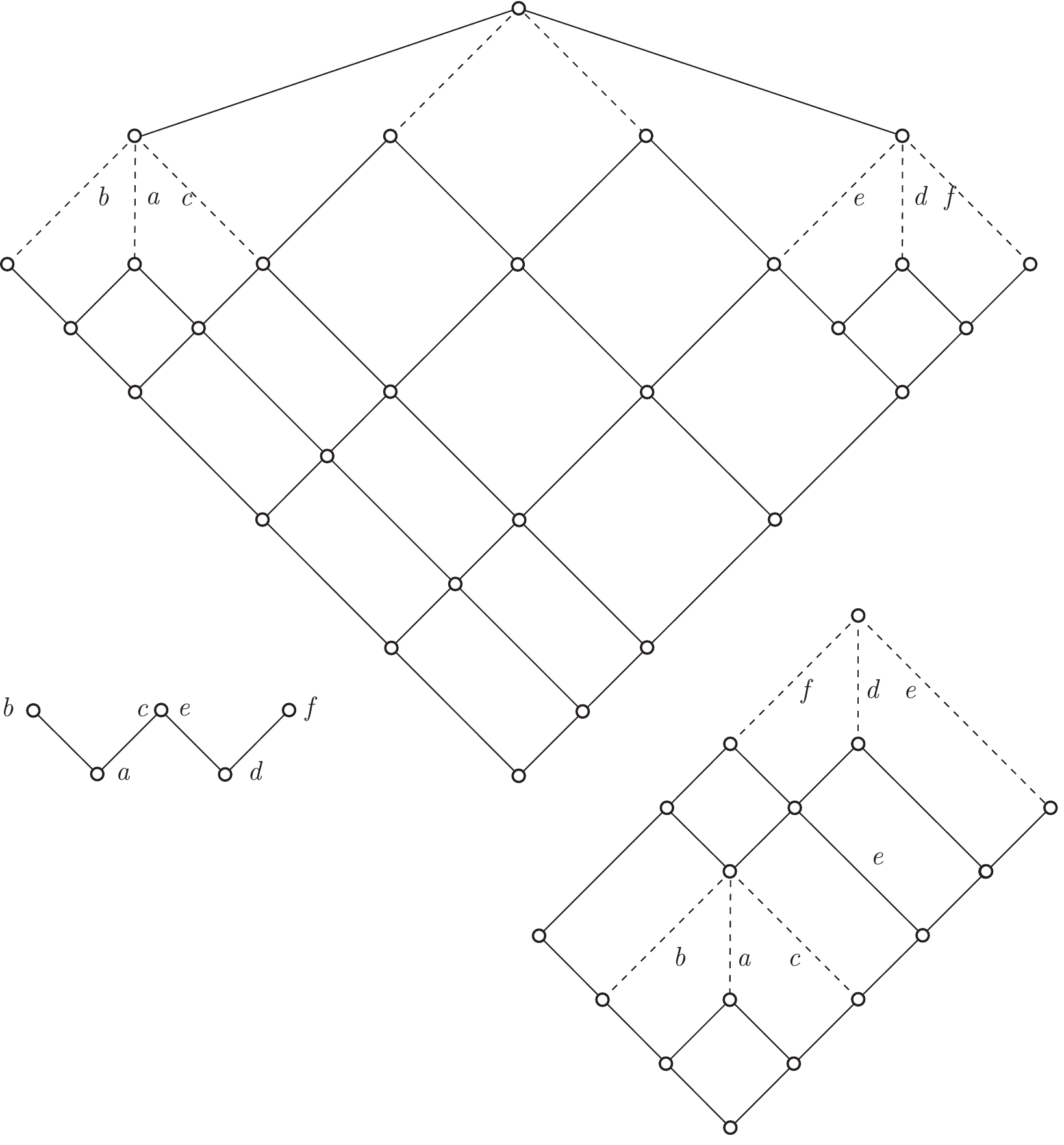}}
\caption{The $W$-relation, two variants}\label{F:W} 
\end{figure}

In this variant, we get $c = e$ by applying the Equality Lemma 
as in the top diagram of Figure~\ref{F:W}.
The dotted edges of the top diagram in Figure~\ref{F:W} are covering 
and they are all distinct except maybe for the two internal edges
at the unit element. Figure~\ref{F:3C} shows the smaller variant of the top diagram
in which these two edges are equal.
The smaller variant will be used in the subsequent diagrams;
the full diagram is too big to draw.

\subsubsection{The W-relation, Version 2}\label{S:Wsecond}

As shown in Figure~\ref{F:W}.

\subsection{The $3$C-relation}\label{S:3C}
The \emph{$3$-Crown-relation}, $\threeC$, in $K$ 
is also a \text{$6$-ary} relation on the set of colors of $K$. 
The relation $V(a,b,c,d,e,f)$ holds in $K$, if 
\[
   V(a,b,c),\  V(d,c,f),\  V(e,b,f)
\] 
hold in $K$, see the diagram on the left of Figure~\ref{F:3C}.

\subsubsection{The $3$C-relation, Version 1}\label{S:3Cv1}

\begin{figure}[b!] 
\centerline{\includegraphics{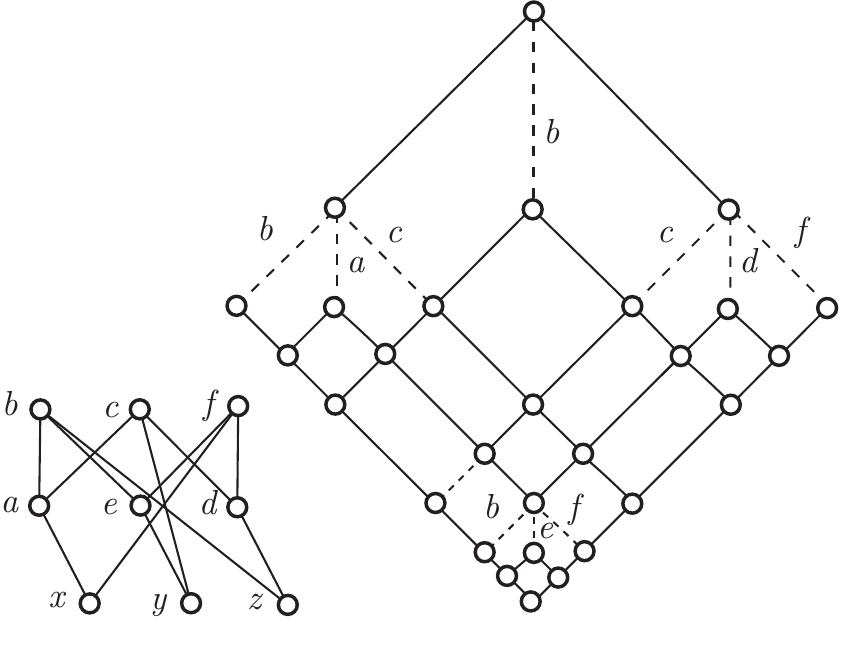}}
\caption{The $\threeC$-relation, Version 1, small} 
\label{F:3C}
\end{figure}

$ V(a,b,c), V(d,c,f)$ form a W-relation, 
so we get the top half of the diagram of the lattice in Figure~\ref{F:3C} (small variant).
 
In addition, $V(e,b,f)$ holds in $K$.
So there is a covering square colored by $b$ and~$f$, 
and we get $V(e,b,f)$  by having a fork (see Section~\ref{S:Forks}) in the covering square.
There is only one covering square in $K$ colored by~$b,f$, 
at the bottom of the diagram.
So we have a fork in there, the middle edge colored by $e$. 
In~larger examples, there may be many covering squares in $K$ 
colored by~$b,f$, we pick one.

\subsubsection{The $3$C-relation, Version 2}\label{S:3CV2}

\begin{figure}[htb]
\centerline{\includegraphics{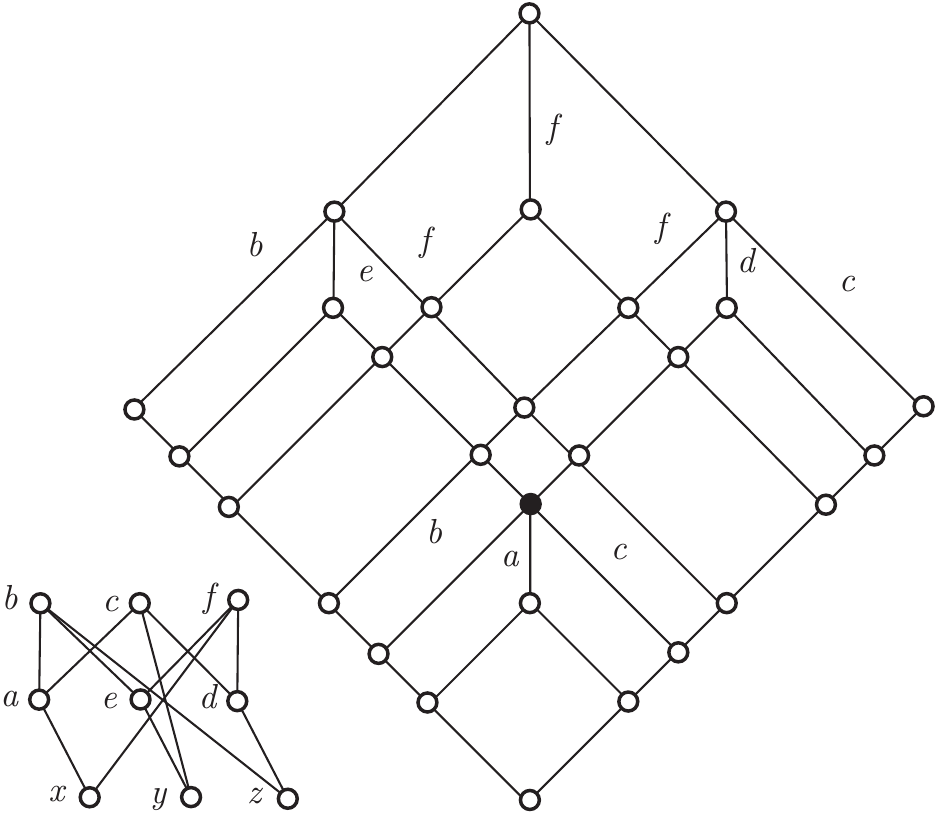}}
\caption{The $\threeC$-relation, Version 2}
\label{F:gotit}
\end{figure}

We need a cover-preserving SPS extension of the bottom lattice of Figure~\ref{F:W}
establishing the W-relation, Version 2, with a cover-preserving square colored by $b$ and $f$,
see Figure~\ref{F:gotit}. The elements of this square are black-filled.
We then insert a fork  (see Section~\ref{S:Forks} into this square.

\section{Proving the 3P3C Theorem}\label{S:theorem}

In my joint paper~\cite{GKn09} with E.~Knapp, we proved that 
an SPS lattice $K$ has a congruence-preserving rectangular extension $\hat K$.
So to verify the 3P3C Theorem for an SPS lattice $K$, 
it is sufficient to verify it for $\hat K$.
Therefore, we can assume that $K$ is an SR lattice.

We want to prove that the ordered set $R_3$ 
has no cover-preserving embedding into $\Ji{\Con K}$.
The relation $\threeC(a,b,c,d,e,f)$ holds in $K$, 
as well, as the three relations
\[
   V(x,a,f),\  V(y,e,c),\  V(z,b,d).
\]
   
\subsection{Version 1}\label{S:V1}

There is an edge $y$ that was obtained by inserting a fork 
in a covering square colored by $c$ and $e$.
But there is no such covering square possible because the $c$ colored edges and 
the $e$ colored edges run parallel.
This completes the proof for Version 1.

\subsection{Version 2}\label{S:V2}

In Figures~\ref{F:gotit}, 
let $I$ denote the ideal generated by the black-filled element.
All the nontrivial congruence classes of the color $a$ are in $I$. 
All the nontrivial congruence classes of the color $f$ are outside of $I$. 
But $x \leq a$ and $x \leq f$, a contradiction. This completes the proof for Version 2.

We have now completed the proof of the 3P3C Theorem.

\end{document}